# Vehicle Routing Problem for Urban and Grey Zones Considering Heterogeneous Fleets and Carpooling

M. Keshvarinia[1], S. M. J. Mirzapour Al-e-hashem[1], R. Shahin[2], M. Kiaghadi[1]

[1]Department of Industrial Engineering, Amirkabir University of Technology, Tehran, Iran

[2] 2500 Chem. de Polytechnique, Montréal, QC H3T 1J4, Canada

**Abstract**

The conveyance of employees holds paramount significance for expansive corporations. Employees typically commute to their workplaces either via personal vehicles or through public transit. In this research endeavor, our role is that of a third-party entity entrusted with orchestrating the transportation of employees whose place of employment is situated within the grey zone. This zone exclusively permits the ingress of electric/hybrid vehicles and buses. We advocate for employees to adopt carpooling and furnish bus services for those who abstain from it. The primary objective of this research is to curtail the quantity of vehicles leased by the third-party entity, promote carpooling among employees, amplify employee contentment, and mitigate environmental degradation stemming from vehicular gasoline consumption. To decipher the model delineated in this study, the epsilon constraint method is proffered for petite-scale instances, while NSGA-II is introduced as a potent meta-heuristic technique tailored for large-scale scenarios. Computational trials corroborate that the models posited can be efficaciously harnessed by enterprises to pare down transportation expenditures.

**Keywords:** Vehicle routing problem, Carpooling, Grey zone, Fleet mixture, NSGA-II algorithm, Epsilon constraint method

## 1. Introduction

The escalating vehicular utilization in urban areas has precipitated a myriad of challenges, encompassing air and noise contamination as well as traffic congestion. In response to these quandaries, Industry 4.0 has emerged, proffering innovative methodologies applicable to urban governance, one of which is the carpooling system. Carpooling entails the shared use of a singular vehicle by multiple employees, facilitating consistent vehicular access without the necessity of individual ownership. Such a system promotes the congregation of numerous individuals within a single vehicle, thereby enhancing traffic fluidity (Mokhtarzadeh et al. 2020). The intent of this research is to devise and actualize an apt model for carpooling, whilst also integrating a bus service for the workforce of an enterprise situated within the Traffic Restriction Zone (TRZ). The TRZ is a precinct exclusively accessible to electric, hybrid vehicles, and buses. Entry into this zone for other vehicular categories is contingent upon specific stipulations, varying across distinct temporal intervals and weekdays. The employee populace is bifurcated into two cohorts: 1. those opting for carpooling; 2. those abstaining from carpooling. The former cohort is mandated to station their vehicles in designated parking facilities prior to TRZ ingress, subsequently utilizing electric or hybrid modes of transport for their workplace commute. Conversely, the latter group can gravitate towards bus terminals, availing a bus service provisioned by a third-party entity to reach their workplace.

In numerous locales, seamless access to public transit remains elusive, compelling individuals to traverse extended distances to alight at a bus or metro terminal. Consequently, corporate entities have embarked on strategic initiatives to facilitate their workforce's utilization of public transportation. Carpooling emerges as a salient strategy, curtailing the reliance on private vehicles for daily commutes. It is imperative to underscore that carpooling systems are optimally tailored for the consistent travel needs of a collective. A conspicuous lacuna in current solutions pertains to the strategic positioning of employees within the TRZ. Hence, it becomes indispensable to address this conundrum to alleviate employee-related challenges. In this scholarly endeavor, we scrutinize an avant-garde mathematical framework aimed at orchestrating the conveyance of a substantial company's employees situated within a designated grey zone.

**2. Literature Review**

## 2.1. Vehicle Routing Problem (VRP)

Over four decades have elapsed since the seminal introduction of the Vehicle Routing Problem (VRP) (Dantzig and Ramser, 1959). Within their scholarly work, they not only elucidated a tangible application of the problem but also pioneered the inaugural mathematical model for its resolution. A few years later, Clarke and Wright (1964) unveiled a superior heuristic methodology, enhancing the approach delineated by (Dantzig and Ramser, 1959). Stemming from these foundational papers, a plethora of both exact and heuristic methodologies have been promulgated, aiming to ascertain both optimal and approximate solutions for diverse iterations of VRPs. The academic realm boasts an array of literature reviews about VRP. Laporte et al. (1986) proffered an exhaustive exploration singularly focused on precise methodologies for VRP resolution, coupled with an in-depth examination of researches undertaken until the terminus of the 1980s. Numerous scholarly endeavors have delved into both exact and heuristic methodologies, as evidenced by works from (Christofides et al., 1981; Toth and Vigo, 1997; Fisher, 1995; Golden et al., 1998; Laporte, 1992; Magnanti, 1981; Toth and Vigo, 2002). Shahin et al. (2016) introduced a linear programming approach, taking into account the correlation amongst right-hand side parameters. It is also noteworthy to mention the burgeoning interest in the Green VRP, which has captivated the academic community in recent years. In the following, we review recent development in the literature.

In Zachariadis et al., (2015), the main contribution of the paper is the development of a load-dependent VRP model, which considers the weight of the cargo in the routing plans. The objective is to minimize the product of the distance travelled and the gross weight carried, which is closely linked to the energy requirements of the vehicle fleet. The results show that the proposed approach outperforms existing methods in terms of solution quality and computational efficiency. In Tirado and Hvattum (2017), the authors present a novel approach to improve solutions for dynamic and stochastic maritime VRP by incorporating stochastic information through local search variations. The proposed heuristics by appending three local search methods to create initial solutions and improve them. In Moumou et al., (2017), the authors propose an updating Kohonen map algorithm to solve the VRP. The proposed algorithm uses unsupervised competitive neural network concepts and includes some traditional and new problem-solving techniques. The methodology used in the paper includes generating random data and executing the algorithm on a computer system. The authors conclude that the proposed algorithm outperforms existing methods in terms of solution quality and computational time. Artificial Intelligence (AI) greatly improves the solutions to the Vehicle Routing Problem (VRP) by optimizing and automating decision-making processes. Machine learning models have the capability to examine past data in order to forecast future demand and client preferences, facilitating more precise planning of routes and schedules. In addition, the ability to process data in real-time enables machine learning to adjust routes immediately in response to traffic conditions, weather fluctuations, and unforeseen delays, so ensuring the most efficient path. An interesting application, as examined in the study of Kiaghadi et al., (2024) involves utilizing machine learning to forecast energy efficiency. This approach can also be applied in a similar manner to predict vehicle fuel efficiency and maintenance requirements in vehicle routing problem (VRP) solutions. The work of Kiaghadi and Hoseinpour (2023), mirrors the utilization of machine learning in vehicle routing problems (VRP) to optimize delivery schedules and routes by employing prescriptive analytics to enhance performance, using several different machine learning models and selecting the best performing model. This highlights the significant impact of machine learning in addressing various logistical challenges. These studies offer useful knowledge on how predictive and prescriptive analytics can improve operational efficiency, which can directly enhance VRP outcomes by optimizing routes and allocating resources intelligently.

In Karak and Abdelghany (2019), the authors develop the hybrid Clarke and Wright heuristic (HCWH), for the hybrid VRP (HVDRP) in pick-up and delivery services. In Beraldi et al., (2021), the authors propose a novel mathematical model for the multi-vehicle VRP with mandatory and not mandatory customers that have to be selected according to given profits. The authors also designed a meta-heuristic algorithm based on a pre-selection process of a sub-set of feasible requests in order to solve small and real-size instances of the problem in a short computational time. The results obtained are validated by the stakeholders as satisfactory.

**2.2. Carpooling**

Amidst the escalating frequency of trips and the prevalent reliance on private vehicles for quotidian activities, the consequent intensification of traffic congestion has illuminated a salient realization for traffic administrators: the mere augmentation of supply and amenities cannot perpetually be regarded as the sole avenue for enhancing the efficacy of diverse facilities and systems. Consequently, strategies centered on demand management have been proffered to realize the anticipated objectives of urban traffic management. The foundational tenets of these strategies can be articulated through the ensuing three overarching paradigms.

- Curtailing the number of trips and excising specific segments thereof.
- Strategically allocating journey timings throughout varied intervals of the day.
- Transitioning the modal preference from individualized transportation towards collective conveyance methods.

Carpooling exerts a beneficial influence on the diminution of travel volume, fuel utilization, environmental contaminants, while concurrently augmenting travel velocity and safety, among other advantages. Asghari et al. (2022) scrutinized the behavior of carpooling within the context of bus routing dilemmas, discerning a resultant decrease in air contamination and greenhouse gas emissions. Alumur and Kara (2008) elucidated in their research that the occupancy quotient of private automobiles (denoting the number of commuters per vehicle for each journey) was notably subdued. The median vehicular occupancy across Europe fluctuates between 1.1 and 1.8 passengers. Santos et al. (2011) corroborated that within the United States, the car occupancy ratio mirrors its European counterpart. Furuhata et al. (2013) pinpointed specific impediments for service purveyors, encompassing the formulation of enticing frameworks, orchestrating the ingress and egress of passengers, and fostering a sense of trustworthiness amidst travelers in digital platforms. The orchestration of carpooling is facilitated through an array of channels, inclusive of social media platforms, corporate web portals, sophisticated mobile applications, and travel intermediaries (Stiglic et al., 2015). Zhang et al. (2016) posited that an impetus for carpooling engagement lies in the introduction of incentives that render transportation more appealing to individuals. In their exploration, they probed into the inducement of designating a specific lane for vehicles accommodating multiple passengers. Their findings underscored that leveraging such lanes culminates in a reduction of both cost and duration for the entire trip.

Lotfi and Abdelghany (2022) propose a novel methodology for integrating ride matching and vehicle routing for On-Demand Mobility Services (ODMS) with ridesharing and transfer options. The hybrid heuristic approach enables solving medium to large problem instances in near real-time, providing efficient solutions while satisfying real-time execution requirements. The methodology decomposes the original problem into a master problem and smaller sub-problems that can be solved separately. The results of the experiments conducted to test the methodology demonstrate its effectiveness in solving the problem of ride matching and vehicle routing for ODMS, with the methodology outperforming two other ridesharing frameworks in literature. Guo et al. (2022) study on the optimization of vehicle routing in intercity ride-sharing services, with the aim of maximizing the total profit of the system. The authors propose a variable neighborhood search algorithm with effective neighborhood operators and local search strategies to solve the vehicle routing problem of intercity ride-sharing. Dai et al. (2022) develop an optimization model for autonomous taxi ride-sharing that takes into account energy consumption. The authors optimized the taxi type, taxi path, and depot location to minimize energy consumption while still providing efficient ride-sharing services. They found that the system was effective in reducing energy consumption compared to regular taxi service but underperformed compared to private vehicle travel in terms of fuel consumption. Wang et al. (2023) study a novel joint decision framework that combines deep reinforcement learning with integer-linear programming to optimize ride-sharing with passenger transfer. The proposed framework addresses the challenges of request dispatching, transfer scheduling, and vehicle rebalancing in a ride-sharing model with passenger transfer. The authors conducted experiments on real-world datasets to evaluate the performance of the proposed framework and compared it with existing methods. The results show that the proposed framework outperforms existing methods in terms of efficiency, effectiveness, and scalability. Pouls et al. (2022) propose a novel approach to repositioning idle vehicles in ride-sharing systems, which uses adaptive forecast-driven algorithms to improve vehicle utilization, customer waiting times, and request rejection rates. The approach is based on a combination of demand forecasting, optimization, and simulation techniques, and it adapts to changes in trip demand and vehicle supply over time. The results of the simulation studies on real-world datasets from Hamburg, New York City, Manhattan, and Chengdu show that the algorithm performs well with perfect and naïve demand forecasts and achieves significant improvements over a simple historical average and classical machine learning algorithms such as gradient boosting.

**2.3. Mixture of the transport fleet**

In contemporary times, the adoption of electric and hybrid vehicles within the transportation domain has emerged as a propitious alternative to mitigate environmental degradation and curtail the sector's reliance on oil and petroleum derivatives (Chan et al., 2009). The allure of electric and hybrid vehicles has witnessed an ascendant trajectory, with market penetration reaching 29% in Norway, 6% in the Netherlands, and 1.5% in nations such as China, France, and the United Kingdom. Nonetheless, despite a decade since their foray into the expansive automobile market, these vehicles grapple with formidable challenges, including elevated price points, limited travel ranges, protracted charging durations, and a paucity of charging infrastructure. The global financial burden attributed to air contamination engendered by fuel-driven vehicles is approximated at a staggering $3 trillion annually. Mancini (2017) introduces a model tailored for the conventional VRP, leveraging hybrid

vehicles powered by an amalgamation of batteries and fossil fuels, and employs the largest neighborhood search algorithm for model resolution. Rezgui et al. (2019) posit that an electric vehicular fleet presents a viable substitute for goods conveyance. In this context, diverse charging locales are contemplated, and the fleet is strategically allocated to these stations. The overarching objective of their scholarly contribution is to conceptualize and ascertain the optimal routing paradigm for a contingent of modular electric vehicles, ensuring goods delivery at the most economical expenditure.

Schenekemberg et al. (2022) develop a new model and solution approach. The model considers a mixed fleet of vehicles with different capacities and user requirements and aims to optimize vehicle routes while meeting service constraints. The solution approach combines branch-and-cut (B&C) with biased random-key genetic algorithm with Q-learning (BRKGA-QL) in a shared parallelized framework. Armbrusta et al. (2022) propose a multi-objective optimization model in rural regions, which takes into account various constraints such as break times of drivers, variable time window sizes, and wheelchair transportation. The authors also analyze two different operator modes on large real-world datasets and evaluate the dynamics of optimal sharing rates and user convenience. Souza et al. (2022) study a bi-objective optimization model and a two-stage heuristic for the heterogeneous, based on a real-world application of patient transportation. The methodology used in the paper includes problem formulation, solution method, and dynamism and trade-off analyses.

Mortazavi et al. (2022) develop a new method for constructing feasible solutions using adaptive logit models. The proposed method uses a stepwise greedy approach to create the sequence of passengers' pick-up/drop-off nodes for every available vehicle, and the process of selecting nodes when constructing routes is dictated by three strategies. The probability of success in constructing a feasible solution is estimated using a binary logit model. Mortazavi et al. (2022) propose a linearly decreasing deterministic annealing algorithm for the multi-vehicle version of the problem. The algorithm guarantees a feasible solution and utilizes advanced local search operators to accelerate the search for optimal solutions. The paper also presents a comprehensive comparison approach to evaluate the efficiency of different solution methods, taking into account aspects such as convergence rate, number of failures, and potential impacts of stochasticity.

## 2.4. Electric vehicles

Electric vehicles are automobiles equipped with batteries and electric motors. These vehicles harness rechargeable batteries for energy and utilize the electric motor for propulsion. A salient merit of electric vehicles lies in their superior energy efficiency. Furthermore, they are devoid of pollutant emissions, playing a pivotal role in markedly diminishing greenhouse gas emissions within the transportation realm. The operational expenditure of electric vehicles is inferior to that of combustion-engine vehicles, a disparity attributable to the differential taxation of electricity and fossil fuels. While electric vehicles are emblematic of environmental stewardship, they are not without challenges. These encompass limited battery capacity coupled with elevated costs. The paucity of charging infrastructures also poses a contemporary impediment. As the adoption of electric vehicles burgeons, the predilection for fuel-based vehicles wanes, thereby potentially serving as a catalyst in the reduction of greenhouse gas emissions (Aziz et al., 2016).

Johnsen and Meisel (2022) present a problem variant that allows for interrelated trips, offering a more holistic approach to meeting the travel demands of users. The authors propose a mixed-integer programming (MIP) formulation of the problem and a heuristic solution representation scheme that can accommodate trip interdependencies. They conduct a computational study to evaluate the performance of the proposed heuristic and compare it with existing methods.

Su et al. (2023) developed a novel metaheuristic algorithm, called Deterministic Annealing (DA), for solving the Electric Autonomous VRP. The DA algorithm is designed to minimize total travel time and excess user ride time while accommodating all customer requests and considering the limited range of electric autonomous vehicles. The methodology involves a fragment-based representation of paths, a linear time complexity route evaluation scheme, and a deterministic annealing local search algorithm. Limmer (2023) propose a bilevel large neighborhood search approach for the electric autonomous variant of the problem, which takes into account battery constraints and the option to recharge vehicles at different charging stations. The results of the experiments show that the proposed approach achieves better results than existing methods on common benchmark instances.

**2.5. Hybrid vehicles**

Hybrid vehicles harness dual, distinct energy sources for propulsion. This entails the deployment of an electric motor in tandem with a conventional fuel-driven engine for power generation. In nations such as Iran and Argentina, vehicles have been conceptualized that utilize a gas engine alongside their traditional fuel-based engine. Nevertheless, such vehicles fall outside the purview of this research. The inaugural hybrid vehicle, a fusion of gasoline and electric mechanisms, materialized in 1990. The commercial advent of hybrid vehicles was marked by Toyota in 1997 with the launch of the Toyota Prius. The dawn of the 21st century witnessed the introduction of electric-hybrid vehicles, which could be powered by domestic electricity.

The hybrid vehicles under examination in this scholarly endeavor possess the capability to propel a vehicle utilizing an energy storage mechanism boasting a power output of four kilowatt-hours or greater. Additionally, these vehicles are equipped with the capacity to recharge their batteries via an external source through a plug interface and can traverse a minimum of ten miles exclusively in electric mode without the consumption of fossil fuels. This particular breed of hybrid vehicle amalgamates the attributes of traditional hybrid vehicles (comprising both fuel-driven and electric motors) and pure electric vehicles. The financial outlay associated with the electricity consumption for this vehicle category is projected to be approximately 1.4 times that of the gasoline expended (Weiss et al., 2019).

In Mancini (2017), the main contribution of the paper is twofold. First, the authors introduce a new and more realistic extension of the Green VRP, the Hybrid VRP, in which vehicles may use both electric and traditional propulsion. Second, they propose an efficient Large Neighborhood Search Meta-heuristic to solve the problem. In Yu (2017), the main contribution of the paper is the proposed mathematical model and simulated annealing with a restart strategy to solve the Hybrid Vehicle Routing Problem (HVRP) with the objective of minimizing the total cost of travel by driving Plug-in Hybrid Electric Vehicles (PHEVs). The model considers the utilization of electric and fuel power, as well as the vehicle's battery capacity and charging time. The simulated annealing with a restart strategy is used to find the optimal solution for the HVRP. The methodology used in the paper includes a literature review, problem definition, mathematical modeling, simulated annealing with a restart strategy, and computational experiments to test the proposed methods. The results show that the proposed methods outperform existing methods in terms of solution quality and computational time. In the following, we propose a literature review summary in Table 1.

**Table 1: Summary of the literature; D (Dynamic); S (Static); Hete (Heterogeneous); Homo (Homogeneous); E (Electric); H (Hybrid); T (Traditional).**

| Ref. | Data | | Environment | | Vehicle | | | | |
|---|---|---|---|---|---|---|---|---|---|
| | Stochastic | Deterministic | D | S | Hete | Homo | E | H | T |
| Lotfi and Abdelghany (2021) | x | | x | | x | | | | x |
| Guo et al., (2022) | | x | x | | x | | | | x |
| Dai et al., (2022) | x | | x | | x | | | | x |
| Wang et al., (2023) | x | | x | | x | | | | x |
| Pouls et al., (2022) | x | | x | | x | | x | | |
| Schenekemberg et al., (2022) | x | | x | | x | | | | x |
| Armbrust et al., (2022) | | x | x | | x | | | | x |
| Souza et al., (2022) | | x | x | | x | | | | x |
| Mortazavi et al., (2022b) | | x | x | | x | | | | x |
| Mortazavi et al., (2022a) | | x | x | | x | | | | x |
| Johnsen and Meisel (2022) | | x | x | | x | | | | x |

| | | | | | | | | |
|---|---|---|---|---|---|---|---|---|
| Su et al., (2023) | x | | x | | x | | x | x | |
| Limmer (2023) | | x | x | | x | | x | | |
| Mancini (2017) | | x | | X | x | | x | x | |
| Yu et al., (2017) | x | | x | | x | | x | x | |
| Tirado and Hvattum (2017) | x | | x | | x | | | | x |
| Beraldi et al., (2021) | x | | x | | x | | | | x |
| Karak and Abdelghany (2019) | | x | x | | x | | | | x |
| **Our study** | | x | | X | x | | x | x | x |

Despite the extensive body of research dedicated to carpooling and the enhancement of employee satisfaction coupled with cost reduction, certain lacunae persist within the extant literature. This scholarly endeavor endeavors to bridge some of these identified voids. In pursuit of this, we amalgamate notions such as carpooling, the grey zone, and advocate for the deployment of electric and hybrid fleets. A model is proffered that calculates the expenditure associated with utilizing each fleet type leased from a tertiary entity and prognosticates the minimal quantity of rental vehicles requisite for employee conveyance. Moreover, by integrating incentives for carpooling, we aim to galvanize employees towards its adoption. Furthermore, through judicious strategizing and the incorporation of electric and hybrid fleets within the TRZ, we anticipate an augmentation in employee contentment and a concomitant diminution in environmental degradation

To this end, the contribution of this research are as follows:

1. Considering the grey zone in a vehicle routing problem.
2. Considering electric and hybrid fleets in grey zone.
3. Integrating multiple VRP in different zones.
4. Using heterogeneous fleet for employees
5. Presenting a multi-objective mathematical model and comparing two proposed solution methods.

In the remainder of the paper, we describe the problem and provide mathematical modeling in Section 2. In Section 3, we provide methodology. Results are presented in Section 4. We conclude in Section 5 and provide future research directions.

## 3. Problem Description

Globally, a significant portion of the populace relies on private automobiles for commuting to work. However, in the contemporary era, given the burgeoning population and ensuing traffic conundrums, carpooling has garnered the attention of numerous corporate entities. The challenges currently plaguing metropolises and densely populated cities predominantly revolve around atmospheric contamination and traffic congestion. These pressing issues and their ramifications have necessitated preemptive strategizing to facilitate more fluid vehicular movement. A salient strategy in this context is the implementation of specific traffic regulations. Pursuant to this regulatory framework, combustion-engine vehicles will be proscribed from accessing the TRZ, with entry being exclusive to electric and hybrid vehicles, as well as buses. As delineated in the extant literature, the quandaries associated with urban employee transportation encompass not only traffic but also pollution, attributable to the prodigious fuel consumption of traditional vehicles and the concomitant restrictions on TRZ access for certain vehicles. Furthermore, the conveyance of employees represents a financial burden for corporations, necessitating judicious management and optimization. In light of these challenges, this research endeavors to proffer solutions that ameliorate the prevailing circumstances and address the aforementioned issues.

To mitigate atmospheric contamination, champion carpooling amongst employees, and augment their contentment, this research delineates the vehicular transportation route into two distinct segments: urban and grey zones. The urban zone is demarcated as an area exclusively accessible to buses and vehicles engaged in carpooling. Conversely, the grey zones, characterized by restrictions on the ingress of combustion-engine vehicles — a measure instituted in response to strategic traffic plans aimed at curtailing pollution and congestion — permit only buses and electric or hybrid fleets.

Regulatory stipulations within this zone dictate that individuals station their vehicles at designated TRZ entry points within parking facilities, subsequently continuing their journey via electric or hybrid vehicles. Given the inherent battery constraints of electric vehicles, rendering them suitable for shorter traverses, hybrid vehicles are preferred for extended commutes within the grey zones. Given the substantial number of individuals grappling with challenges related to TRZ access on a routine basis, corporations situated therein have elected to ameliorate these issues through collaboration with tertiary entities specializing in transportation management. Within the purview of this research, it is postulated that an enterprise has delegated the conveyance of its workforce to such a third-party entity, with the overarching objectives of curtailing transportation expenditures, diminishing environmental degradation, and bolstering employee contentment. This tertiary firm proffers incentives to those opting for carpooling, categorizing them based on geographical proximities. For those abstaining from carpooling, buses are provisioned at diverse urban locales. Employees possessing vehicles, in a display of voluntarism, accommodate co-workers during their commute, thereby economizing on their expenses. In recognition of this, the third-party entity extends carpooling inducements. An illustrative incentive might be the waiver of parking fees in the vicinity of TRZ ingress points. Employees converge at stations proximate to these parking facilities, subsequently accessing the TRZ via electric and hybrid vehicles leased by the third-party firm. To elucidate further, this research contemplates two distinct segments:

1. Employees either go to the TRZ entry points through the carpooling or they will go to the nearest bus station and get on the rented buses by a third-party company. It should be noted that in this study, the location of TRZ entry points and bus stations is clear and people do not get off the vehicles in the middle of the route. The objectives of this study can be summarized as follows:

• Use the maximum capacity of each vehicle along the route to the TRZ entry points.

• Determine that each employee has to be assigned to which vehicle, and to determine the pick-up of which employee should be at which station and through which route arrive at the drop-off station in a way to reduce the waiting time for employees and the greenhouse gas emissions.

• Reduce walking time for employees who use buses to get to the TRZ entry points.

2. The second part is related to the grey zone. In this section, employees who have come from home through the carpooling to the entry points of TRZ, depending on the distance from the station to their office, get on a/an electric/hybrid vehicle rented by a third-party company and get to work. The objectives of this study are as follows:

   • Optimal fleet allocation according to the distance of people from the stations to the company

   • Minimize the number of vehicles

   • Reduce fixed fleet costs

   • Minimize the waiting time for grey zone vehicles at stations.

   • Reduce environmental pollution by allocating the grey zone fleet to the stations in such a way that they use their maximum capacity and reach the company with the least energy consumption.

It should be noted that in this study, the model selects to which station employees should be assigned to, as well as which electric and hybrid fleet pick them up and drop them off to the company.

In alignment with the stipulated objectives, we employ mathematical frameworks and optimization techniques with the intent to curtail greenhouse gas emissions, contingent upon the fuel consumption metrics of individual vehicles. Concurrently, we aim to diminish the waiting durations experienced by customers at diverse junctures, thereby enhancing their overall satisfaction. Furthermore, we endeavor to truncate passenger travel durations by selecting the most direct routes between points, and to pare down fixed expenditures, given our mandate overseeing the electric and hybrid fleet pertinent to the grey zone. This research delineates the specific routes wherein electric and hybrid fleets are deployed within this zone. To this end, the research objectives are categorized into three distinct facets: Minimization of emissions from fuel-driven fleets, reduction of the aggregate cost associated with the leased fleet, and curtailing employee waiting intervals, which invariably augments their satisfaction.

This investigation accentuates the VRP paradigm, incorporating carpooling facilitated by buses, electric, and hybrid fleets, under the ensuing assumptions: The urban zone's fleet comprises fuel-driven vehicles, whereas the grey zone predominantly utilizes electric and hybrid vehicles. Buses serve as a ubiquitous mode of transport across both zones. Within the urban precinct, there exists a cohort keen on vehicle sharing or carpooling. Electric fleets cater to shorter commutes, while hybrid fleets are designated for more extended routes within the grey zone. Carpooling participants commence their journeys from their residences (a door-to-door model), progressing to the TRZ entry points, with no en-route disembarkations. These carpooling individuals also avail themselves of complimentary parking amenities. The study contemplates unidirectional journeys, focusing solely on the employees' commute to their workplace, excluding the return leg. The transportation requisites of all employees fall under the purview of the third-party entity. All parameters are unequivocally defined. Buses, stationed at specific urban locales, ferry employees to their workplaces. Electric and hybrid fleets, in conjunction with buses, are leased by a third-party entity, and their usage is exclusive to the company's employees. During their commute, employees remain aboard a singular vehicle, ensuring no transitions between zones. Only those engaged in carpooling undergo a vehicular switch at the TRZ entry juncture. For better understanding, we have provided the problem description in Figure 1.

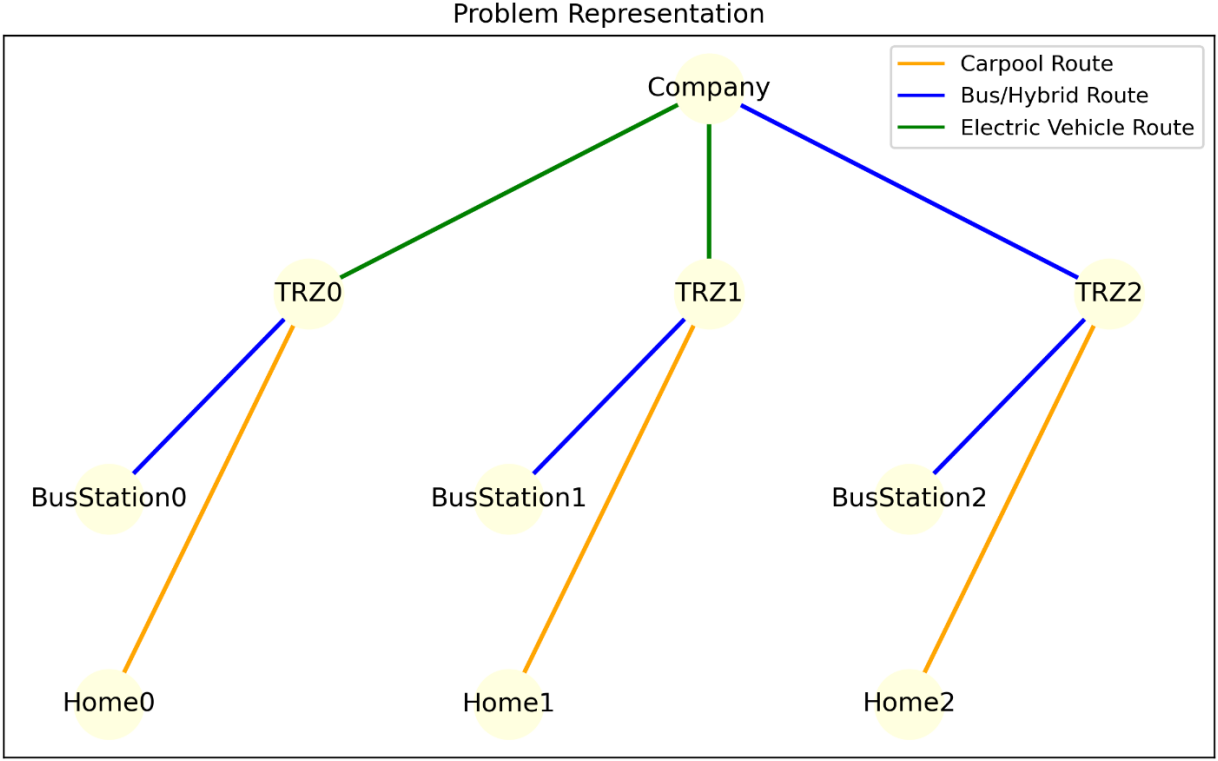

**Figure 1: Problem description**

## 4. Mathematical Formulation

We model the problem on a directed graph where each customer or stations are considered to be a node and there are arcs connecting the nodes. In the following we present the notation.

**Sets**:
- $N$    All nodes of employees ($W$), bus stations ($V$) entry points of TRZ ($S$) and company location ($M$).
- $C$    Set of employee nodes that carpooling is possible (have a personal car).
- $V$    Set of bus stations nodes for v ∈ V, where 0 is the origin and M is the company (destination)
- $S$    Set of nodes that electric and hybrid vehicles stop. for $s \in S$ and $0^1$ is the starting node and M is the company (destination node).
- $K$    Set of all vehicles for k ∈ K
- $B$    Set of buses for b ∈ B.
- $G$    Set of fuel-based vehicles for g ∈ G.
- $E$    Set of electric vehicles for e ∈ E.
- $H$    Set of hybrid vehicles for h in H.

**Parameters**:
- $t_{ij}^k$    Travel time between $i$ and $j$ for $i, j \in N$ by vehicle $k \in K$ (in minutes).
- $FP^k$    Fixed vehicle cost $k \in K$ for hybrid and electric vehicles and buses.

| | |
|---|---|
| $\theta$ | Satisfaction (as a negative parameter) or dissatisfaction (as a positive parameter) for each kilometer walked by each employee |
| $\Delta$ | Time at which working time starts and employees must be at the company before that. |
| $Q^k$ | Capacity of vehicles $k \in K$ |
| $c^k$ | Variable cost of vehicle driver ($k \in K$). |
| $L$ | A big number and is set to the overall of the system to be able to satisfy the constraints |
| $d_{ij}$ | Distance traveled from $i$ and $j$, for $i.j \in C$ (in minutes). |
| $dg_{ij}$ | The distance that a hybrid vehicle travels with fossil fuels for $i, j \in S$. |
| $\beta$ | Carpooling incentive. |
| $r$ | Battery consumption rate in electric and hybrid vehicles per kilometer. |
| $Qe$ | Battery capacity of each vehicle. |
| $t^w$ | The average time it takes for a person to walk one kilometer (in hours) |
| $\varepsilon$ | dissatisfaction due to vehicles arriving earlier or late. |
| $FC^k$ | Vehicle carbon dioxide emissions $k \in K$ (g / km) |
| $de$ | The maximum distance that electric vehicles can travel. |

**Decision Variables**:

| | |
|---|---|
| $x_{ijk}$ | Equals to one if vehicle $k$ goes from $i$ to $j$. |
| $Y_{ij}$ | Equals to one if a person goes from one node $i \in W$ to another node $j \in V$ to get on the bus |
| $Z_j^k$ | The battery level of electric and hybrid vehicles when it arrives at node $j$ |
| $H_k$ | Equals to one if vehicle $k$ is used in the fleet. |
| $H2_i$ | Equals to one If at least one employee is assigned to one bus station. In the case of entry points in the TRZ, if at least one vehicle is assigned to that station. |
| $t_i^k$ | The time vehicle $k$ arrives at node $i$. |
| $t1^k$ | To compute the rate of early arrival at the company (in minutes). |
| $t2^k$ | To compute the rate of late arrival at the company (in minutes) |
| $u_{ik}$ | Sub-tour elimination from the vehicles' route |
| $xy_{ij}^k$ | To linearize capacity constraints of buses. |
| $H3_i$ | To indicate the number of people reaching each TRZ entry point |
| $sx_{ijs}^k$ | To linearize constraints related to computing the number of people entering each station. |
| $H3x_{ij}^k$ | To linearize capacity constraint in TRZ. |

$$Min \sum_{k \in K \backslash G} \sum_{i,j \in N} c^k d_{ij} x_{ijk} + \sum_{k \in (K/G)} FP^k H_k + \beta \sum_{k \in G}(1 - H_k) \quad (1)$$

$$Min \sum_{i \in W} \sum_{j \in V} \theta t^w d_{ij} Y_{ij} + \sum_{k \in B \cup H \cup E} \varepsilon(t1^k + t2^k) \quad (2)$$

$$Min \sum_{i,j \in W} \sum_{k \in G \cup B} FC^k d_{ij} x_{ijk} + \sum_{i,j \in S \cup M} \sum_{k \in H} FC^k dg_{ij} x_{ijk} \quad (3)$$

$$\sum_{j \in V} x_{0jk} = H_k, \forall k \in B \quad (4)$$

$$\sum_{i \in V} x_{iMk} = H_k, \forall k \in B \quad (5)$$

$$\sum_{\substack{i \in 0 \cup V \\ i \neq p}} x_{ipk} - \sum_{\substack{j \in S \cup V \\ j \neq p}} x_{pjk} = 0, \forall k \in B, p \in V \quad (6)$$

$$u_{jk} - u_{ik} + |V|x_{ijk} \leq |V-1|, \quad \forall i \in V, j \in V, k \in B, i \neq j \tag{7}$$

$$\sum_{k \in B} \sum_{\substack{i \in 0 \cup V \\ i \neq j}} x_{ijk} = H2_i, \forall j \in V \tag{8}$$

$$\sum_{i \in 0 \cup V} \sum_{j \in V} x_{ijk} \sum_{l \in W} Y_{lj} = Q^k H_k, \forall k \in B \tag{9}$$

$$\sum_{j \in V} Y_{ij} + \sum_{k \in G} \sum_{\substack{j \in W \\ i \neq j}} x_{ijk} = 1, \forall i \in W \tag{10}$$

$$\sum_{i \in W} Y_{ij} \leq LH2_j, j \in V \tag{11}$$

$$xy_{ijk} \leq \sum_{l \in W} Y_{lj}, \forall k \in B, i \in 0 \cup V, j \in V \tag{12}$$

$$xy_{ijk} \leq Lx_{ijk}, \forall k \in B, i \in 0 \cup V, j \in V \tag{13}$$

$$xy_{ijk} \geq \sum_{l \in W} Y_{lj} - L(1 - x_{ijk}), \forall k \in B, i \in 0 \cup V, j \in V \tag{14}$$

$$\sum_{i \in 0 \cup V} \sum_{j \in V} xy_{ijk} = Q^k H_k, \forall k \in B \tag{15}$$

$$\sum_{\substack{j \in W \\ j \neq C}} x_{cjk} = H_k, \forall c \in C, k \in G \tag{16}$$

$$\sum_{i \in W} \sum_{j \in S} x_{ijk} = H_k, \forall k \in G \tag{17}$$

$$\sum_{\substack{i \in W \\ i \neq p}} x_{ipk} - \sum_{\substack{j \in S \cup W \\ j \neq p}} x_{pjk} = 0, \forall k \in G, p \in W \backslash C \tag{18}$$

$$u_{jk} - u_{ik} + |C|x_{ijk} \leq |C-1|, \quad \forall i \in W, j \in W, k \in G, i \neq j \tag{19}$$

$$2 \leq \sum_{i,j \in W} x_{ijk} \leq Q^k, \forall k \in G \tag{20}$$

$$t_i^k + t_{ij}^k \leq t_j^k + L(1 - x_{ijk}), \forall i \in V \cup 0, j \in V \cup M, \forall k \in B \tag{21}$$

$$t_i^k + t_{ij}^k \geq t_j^k - L(1 - x_{ijk}), \forall i \in V \cup 0, j \in V \cup M, \forall k \in B \tag{22}$$

$$t_i^k + t_{ij}^k \leq t_j^k + L(1 - x_{ijk}), \forall i \in W, j \in W \cup S, \forall k \in G \tag{23}$$

$$t_i^k + t_{ij}^k \geq t_j^k - L(1 - x_{ijk}), \forall i \in W, j \in W \cup S, \forall k \in G \tag{24}$$

$$\sum_{k \in G} \sum_{i \in W} x_{ijk} \leq LH2_j, \forall j \in S \tag{25}$$

$$\sum_{j \in S \cup M} x_{0^1 jk} = H_k, \forall k \in E \cup H \tag{26}$$

$$\sum_{i \in S} x_{iMk} = H_k, \forall k \in E \cup H \tag{27}$$

$$\sum_{\substack{i \in 0 \cup S \\ i \neq p}} x_{ipk} - \sum_{\substack{j \in S \cup M \\ j \neq p}} x_{pjk} = 0, \forall k \in E \cup H, p \in S \tag{28}$$

$$u_{jk} - u_{ik} + |S|x_{ijk} \leq |S-1|, \quad \forall i \in S, j \in S, k \in E \cup H, i \neq j \tag{29}$$

$$\sum_{k \in E \cup H} \sum_{\substack{i \in 0 \cup S \\ i \neq j}} x_{ijk} = H2_j, \forall j \in S \tag{30}$$

$$\sum_{i \in S \cup 0} \sum_{\substack{j \in S \cup M \\ i \neq j}} x_{ijk} d_{ij} \leq d_e, \forall k \in E \tag{31}$$

$$t_i^k + t_{ij}^k \leq t_j^k + L(1 - x_{ijk}), \forall i \in S \cup 0, j \in S \cup M, \forall k \in H \cup E \tag{32}$$

$$t_i^k + t_{ij}^k \geq t_j^k - L(1 - x_{ijk}), \forall i \in S \cup 0, j \in S \cup M, \forall k \in H \cup E \tag{33}$$

$$t_i^k \geq t_i^{k1} - L\left(1 - \sum_{j \in W} x_{jik1}\right), \forall i \in S, \forall k \in H \cup E, k1 \in G \tag{34}$$

$$t1^k \geq \Delta - t_M^k - L(1 - H_k), \forall k \in B \cup E \cup H \tag{35}$$

$$t2^k \geq t_M^k - \Delta - L(1 - H_k), \forall k \in B \cup E \cup H \tag{36}$$

$$\sum_{k \in G} \sum_{i \in W} \sum_{j \in W \cup S} x_{ijk} \sum_{l \in W} x_{lsk} \leq H3_i, \forall s \in S \tag{37}$$

$$\sum_{i,j \in S \cup M} H3_i x_{ijk} \leq Q^k, \forall k \in H \cup E \tag{38}$$

$$sx_{ijsk} \leq \sum_{l \in W} x_{lsk}, \forall k \in G, i \in W, j \in W \cup S, s \in S \tag{39}$$

$$sx_{ijsk} \leq Lx_{ijk}, \forall k \in G, i \in W, j \in W \cup S, s \in S \tag{40}$$

$$sx_{ijsk} \geq \sum_{l \in W} x_{lsk} - L(1 - x_{ijk}), \forall k \in G, i \in W, j \in W \cup S, s \in S \tag{41}$$

$$\sum_{k \in G} \sum_{i \in W} \sum_{j \in W \cup S} sx_{ijsk} \leq H3_i, \forall s \in S \tag{42}$$

$$H3x_{ijk} \leq H3_i, \forall k \in H \cup E, i \in S, j \in S \cup M \tag{43}$$

$$H3x_{ijk} \leq Lx_{ijk}, \forall k \in H \cup E, i \in S, j \in S \cup M \tag{44}$$

$$H3x_{ijk} \leq H3_i - L(1 - x_{ijk}), \forall k \in H \cup E, i \in S, j \in S \cup M \tag{45}$$

$$\sum_{i,j \in S \cup M} H3x_{ijk} \leq Q^k, \forall k \in H \cup E \tag{46}$$

$$Z_0^k = QeH_k, \forall \, \forall k \in E \cup H \tag{47}$$

$$Z_M^k \geq r(d_{M0} - dg_{M0})H_k, \forall k \in H \cup E \tag{48}$$

$$Z_j^k \leq Z_i^k - r(d_{ij} - dg_{ij}) + L(1 - x_{ij}^k) \, \forall k \in E \cup H, \forall i \in S \cup 0, \forall j \in S \cup M \tag{49}$$

$$x_{ijk}, Y_{ij}, H_k, H2_i \in \{0,1\}. \tag{50}$$

$$Z_j^k, t_i^k, t1^k, t2^k, u_{ik}, xy_{ijk}, H3_i, sx_{ijsk}, H3x_{ijk} \geq 0 \tag{51}$$

The first objective function (1) computes the fixed and variable cost of bus, electric and hybrid fleets drivers according to the distance traveled from $i$ to $j$ and maximizes the amount of incentives given to the drivers of fuel-based vehicles. It is notable that it is assumed that not all of the employees have a vehicle to use to get to the office. Therefore, the third component of the first objective function allocate incentives to those use their vehicles. The second objective function (2) by allocating employees to the nearest bus stop minimizes their walking time. Also, it minimizes the arriving time at work early or late for employees by considering the penalty co-efficient for fleets that do not have a proper schedule to arrive to the company. The third objective function (3) minimizes greenhouse gas emissions from buses, fuel based and hybrid vehicles. Constraints (4) to (8) are for routes between bus stops. Constraints (4) states that each bus must start from node zero. Constraints (5) indicates that any bus used in the transport fleet must go to the company node. Constraints (6) states that if a bus enters a bus stop, it must leave that node. Constraints (7) eliminate the sub-tour from the bus route. Constraints (8) states that each bus stop should be visited at most once. Constraints (9) indicates the capacity of buses. Constraints (10) states that each employee must be assigned to either a bus stop or use carpooling. Constraints (11) indicates whether employees are assigned to a bus stop. Constraints (12) to (14) are related to linearization of $x_{ij}^k \sum_{l \in W} Y_{lj}$. Constraints (9) is replaced by constraint (15). Constraints (16) to (21) are for the vehicles route. Constraints (16) states that each carpooling vehicle must start moving from its own starting node. Constraints (17) states that each vehicle must enter exactly one node of the TRZ entry points. Constraints (18) states that if a vehicle enters a node, it must leave that node. Constraints (19) is for sub-tour eliminations of nodes. Constraints (20) is for the capacity of carpooling vehicles. Constraints (21) to (24) computes the time to reach each node in the urban zone. Constraints (21) and (22) are related to buses and Constraints (23) and (24) are related to carpooling vehicles. Constraints (25) indicates which of the entry points of TRZ is used. Constraints (26) to (30) are related to routes between stations in TRZ. Constraints (26) states that any electric and hybrid fleet must start from the node $0^1$. Constrains (27) states that any vehicle used in the transport system must eventually go to the company. Constraints (28) states that if a vehicle enters a TRZ entry point, it must leave that node. Constraints (29) is related to the sub-tour elimination from the route of the electric and hybrid fleets. Constraints (30) states that each TRZ entry point should be visited at most once. Constraints (31) indicates that the distance traveled by the electric vehicle must be less than the maximum distance at which it can drive (based on the battery charge). Otherwise, a hybrid car is used. Constraints (32) and (33) are computing the travel time between nodes in the TRZ. Constraint (34) states the start time of travel from each of the entry stations of the TRZ, which must be after the passengers' arrival time at the station. Constraints (35) and (36) computes the arriving time (early or late) to the company. Constraints (37) computes the total number of people entering each TRZ entry points. Constraints (38) is related to the capacity of vehicles in the TRZ. Constraints (39) to (42) are for linearizing Constraints (37). Constraints (37) is also replaced by constraints (42). Constraints (43) to (46) are related to the linearization of Constraints (38). Constraints (38) is also replaced by constraint (46). Constraints (47) to (49) apply to the charge management of electric and hybrid vehicles. Constraints (47) states that if a vehicle is used, it will be fully charged at the starting node. Constraints (48) indicate the energy required for electric and hybrid vehicles to return to the stations. Constraints (49) indicate the battery level of the electric and hybrid fleets at the arriving time at node. ($dg_{ij}$ distance from $i$ to $j$ is the distance traveled by gasoline). Constraints (50) and (51) indicate the decision variables.

## 5. Methodology

Within the confines of this manuscript, we employ the epsilon constraint methodology, wherein a singular objective function is earmarked for optimization, whilst the remaining functions are integrated as constraints of the model. An exhaustive elucidation of the Epsilon constraint method is furnished in Online Appendix A. Additionally, to tackle the quandary for expansive scale instances, we resort to the NSGA-II meta-heuristic, the intricacies of which are delineated in Section (3.1). The efficacy of NSGA-II has been substantiated (Rabbani et al., 2021). Deliberations concerning parameter tuning and Taguchi analysis are presented in Online Appendix B. We have also provided the description of instances used for computational results in Appendix C. Finally, we provide the details of NSGA-II algorithm in appendix D.

## 5.2 Solutions representation

To ensure an elevated performance caliber for meta-heuristic algorithms, the solutions to the quandary ought to be architectured in a manner conducive to seamless interfacing with the meta-heuristic algorithm. This representation of solutions should encompass the entirety of potential resolutions to the issue while minimizing infeasible modalities pertinent to the problem. Given that the issue at hand is a derivative of the VRP, the representation of problem solutions is facilitated through the utilization of multiple arrays, analogous to the representation methodology employed for VRP solutions. Pertaining to our specific challenge, five distinct arrays are proffered to depict the solutions.

The dimension of the inaugural array is determined by the summation of the total number of employees and the count of bus stops, diminished by one. This array delineates the specific bus station to which each employee is allocated. For instance, given a scenario with 10 employees and 4 designated bus stop locales, the resultant array would span a length of 13. The constituent values of this array would manifest as a permutation ranging from 1 through 13. Interpreting this array, any employee denoted by a number preceding the numerical representation of a bus station will be affiliated with said station. Thus, the numerals 1 through 10 symbolize the employees, while the integers 11, 12, and 13 epitomize the three successive stations. The positioning of an employee's numeral anterior to a bus station's numeral signifies their allocation to that particular station. Conclusively, at the array's terminus, employees are relegated to the final bus station. An illustrative exemplar is proffered in Figure 3. Within this exemplification, the third bus station remains devoid of any employee allocation, as no numeral precedes 13. Employees denoted by 10, 8, and 1 are affiliated with bus station 2; those represented by 7, 9, and 5 are associated with station 1; and the employees identified by 6, 4, 3, and 2 are allocated to the terminal station.

| 13 | 10 | 8 | 1 | 12 | 7 | 9 | 5 | 11 | 6 | 4 | 3 | 2 |
|----|----|---|---|----|---|---|---|----|---|---|---|---|

Figure 3. An example of assigning employees to bus stations

The size of the next array equals the total employees plus the number of carpooling vehicles. Like before, this array arranges all employees and vehicles in a specific order. This order shows how each vehicle will travel and which employees it will pick up. Each vehicle carries as many employees as it can hold. Employees did not assign a vehicle number in the array do not get a carpool spot. After figuring out who will carpool, these carpoolers are removed from their original bus stops. For example, in Figure 4, every vehicle carries four people, including one driver and three passengers. Here, employee 1 drives alone. Employee 2 carpools with employees 8 and 10. Employee 3 carpools with employees 5, 7, and 9. The rest take the bus to work. .

| 1 | 11 | 2 | 10 | 12 | 3 | 9 | 5 | 13 | 6 | 4 | 7 | 8 |

Figure 4. An example of carpooling array

The dimension of the tertiary array is determined by the summation of the total number of bus stops and the count of buses, diminished by one. This array delineates the specific route undertaken by each bus. Should a bus stop lack any assigned employees, it will be excised from the designated route. Additionally, each bus will accommodate employees up to its maximum seating capacity. For illustrative elucidation, one may refer to Figure 5. Within this representation, there are four bus stations and an equivalent number of buses, rendering the array's length as seven. As per this depiction, buses 2 and 1 remain devoid of designated routes (they remain unutilized). Bus 3 traverses through stations 2 and 1, while bus 4 navigates through stations 3 and 4.

| 6 | 2 | 1 | 7 | 5 | 4 | 3 |

Figure 5. An example of a bus route array

The quartile array bears resemblance to the third, albeit it delineates the route from TRZ entry stations to the corporate establishment. Its dimension is the summation of the count of TRZ entry junctures and the aggregate of electric and hybrid vehicles, diminished by one. It is imperative to highlight that a vehicle can navigate to a novel entry station provided it retains sufficient energy reserves. Moreover, should the span to the company be extensive, surpassing the maximum range achievable by electric vehicles, a hybrid vehicle is deployed. Only those stations with allocated employees are incorporated within the route. The quintile array specifies the TRZ entry junctures for carpooling vehicles. This array spans a length equivalent to the count of carpooling vehicles, and each cell's value within this array signifies the numeral of the designated ingress station.

## 6. Results

**Before we dive into our results, here we define the characteristics of our instances in Table 2.**

**Table 2: Details of the instances**

| Indicator | Sign | Instances | | | | | | | | | | | | |
|---|---|---|---|---|---|---|---|---|---|---|---|---|---|---|
| | | I1 | I2 | I3 | I4 | I5 | I6 | I7 | I8 | I9 | I10 | I11 | I12 | I13 |
| Number of bus stops | V | 2 | 2 | 3 | 3 | 4 | 4 | 5 | 5 | 6 | 6 | 8 | 8 | 8 |
| Number of TRZ entrance stops | S | 1 | 3 | 3 | 4 | 4 | 5 | 5 | 6 | 6 | 6 | 8 | 8 | 8 |
| Number of employees | W | 10 | 15 | 20 | 25 | 30 | 35 | 40 | 45 | 50 | 55 | 60 | 65 | 70 |
| Number of ridesharing employees | C | 3 | 4 | 4 | 6 | 6 | 8 | 8 | 12 | 12 | 12 | 15 | 15 | 15 |
| Number of buses | B | 1 | 1 | 1 | 2 | 2 | 2 | 3 | 3 | 3 | 4 | 4 | 5 | 5 |
| Number of fuel-based vehicles | G | 3 | 4 | 4 | 6 | 6 | 8 | 8 | 12 | 12 | 12 | 15 | 15 | 15 |
| Number of electric vehicles | E | 1 | 2 | 3 | 3 | 4 | 5 | 6 | 6 | 7 | 8 | 8 | 9 | 9 |
| Number of hybrid vehicles | H | 2 | 2 | 4 | 4 | 5 | 6 | 7 | 7 | 7 | 8 | 8 | 9 | 9 |
| Covering | Nj | All | All | All | All | All | All | All | All | All | All | All | All | all |

**In the following, we discuss our results.**

**6.1 Comparing the performance of NSGA-II for small-scale instances**

To authenticate the efficacy of the genetic algorithm, the outcomes derived from it are juxtaposed against those procured via the epsilon constraint method. The experiments were executed on a computer with an Intel(R) Core (TM) i5–3337 U CPU @ 1.80 GHz, 7858MiB RAM, and Ubuntu 18.4 operating system and the code was programmed in MATLAB version 2017a. It is notable that the reported results for the exact method reach optimality (0% gaps). A comparative analysis is conducted on three exemplar problems, as delineated in Tables 2 through 4. The resultant data underscores the capability of the proposed algorithm to discern the optimal value for each objective function through the precise method within an acceptable temporal frame. The principle of dominance has been invoked to appraise the algorithm's performance. As evinced in the tables, the Pareto solutions proffered by the genetic algorithm closely align with those of the epsilon constraint method, attesting to the superior caliber of the proposed algorithm in addressing the problem at hand.

Table 3. Pareto solutions obtained by the algorithm and the exact method for the first example

| Pure Solver | | | NSGA-II | | |
|---|---|---|---|---|---|
| First obj. | Second obj. | Third obj. | First obj. | Second obj. | Third obj. |
| 2395.7 | 110172 | 9060 | 4276.1 | 55124 | 42862 |
| 3148 | 110876 | 7280 | 2658.4 | 75916 | 26063 |
| 5545.5 | 56824 | 42832 | 3625 | 58016 | 50687 |
| 3148 | 110876 | 7280 | 2465 | 82196 | 15000 |
| 4960.8 | 47212 | 52482 | 3626.4 | 68012 | 47722 |
| 5300.8 | 73768 | 36656 | | | |

Table 4. Pareto solutions obtained by the algorithm and the exact method for the second example

| Pure Solver | | | NSGA-II | | |
|---|---|---|---|---|---|
| First obj. | Second obj. | Third obj. | First obj. | Second obj. | Third obj. |
| 3754.2 | 61051 | 12904 | 2976.6 | 44748 | 17792 |
| 3266.2 | 66836 | 13184 | 3468.6 | 37895 | 17452 |
| 3513.6 | 24019 | 32920 | 3744.2 | 56868 | 13072 |
| 3513.6 | 24019 | 32920 | 3054.6 | 52507 | 13868 |
| 3754.2 | 61051 | 12904 | 3542.6 | 46722 | 13588 |
| | | | 3252.2 | 63721 | 13412 |

Table 5. Pareto solutions obtained by the algorithm and the exact method for the third example

| Pure Solver | | | NSGA-II | | |
|---|---|---|---|---|---|
| First obj. | Second obj. | Third obj. | First obj. | Second obj. | Third obj. |
| 4417.5 | 57592 | 36376 | 4452 | 63640 | 34394 |
| 4712.2 | 76142 | 15882 | 4929.5 | 37554 | 28596 |
| 5239.6 | 30830 | 43485 | 4519 | 73276 | 21527 |
| 4680.6 | 33804 | 44883 | 4533.5 | 49342 | 17023 |
| 5332.5 | 44092 | 20114 | 4832.2 | 42962 | 27082 |
| 4712.2 | 76142 | 15882 | 4600.6 | 35166 | 37428 |

**6.2 Performance of NSGA-II algorithm for medium and large scale instances**

Within this segment, the efficacy of the advanced algorithm is assessed post parameter tuning, utilizing the meticulously crafted sample problems. The evaluation of the proposed algorithm's performance hinges on four pivotal metrics: Number of Pareto Solution (NPS), Spacing Metrics (SM), Diversification Metrics (DM), and Mean Ideal Distance (MID). The relevant figures are shown from Figures (6) – (9). In these figures, the x-axis represents the instance number and the y-axis represents the corresponding criterion. For example, the y-axis of Figure (6) is SM. Each experimental scenario underwent quintuple testing via the proposed algorithm. Subsequently, the average values of the aforementioned indicators, amassed over the five iterations for each problem, are documented and presented in Table 4.

Table 6. The average value of the obtained indicators for the algorithm during 5 experiments

| Instances. | MID | NPS | DM | SM |
|---|---|---|---|---|
| | NSGA-II | NSGA-II | NSGA-II | NSGA-II |
| 1 | 0.9 | 15 | 1 | 0.3 |
| 2 | 1 | 10 | 1 | 0.3 |
| 3 | 1 | 7 | 0.9 | 0.03 |
| 4 | 0.9 | 17 | 0.9 | 0.2 |
| 5 | 0.9 | 29 | 1 | 0.3 |
| 6 | 0.9 | 37 | 0.9 | 0.2 |
| 7 | 0.9 | 23 | 1 | 0.1 |
| 8 | 1 | 32 | 1 | 0.2 |
| 9 | 0.9 | 32 | 1 | 0.2 |
| 10 | 1 | 22 | 1 | 0.1 |
| 11 | 0.9 | 34 | 0.9 | 0.2 |
| 12 | 1 | 21 | 0.8 | 0.1 |
| 13 | 0.9 | 27 | 0.8 | 0.07 |

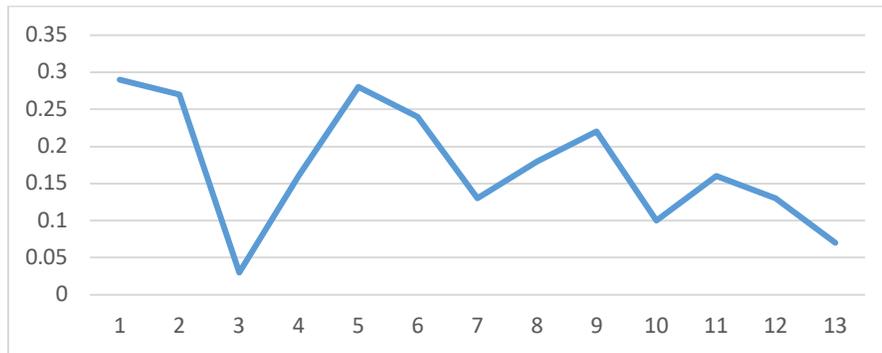

Figure 6. The mean values of the SM

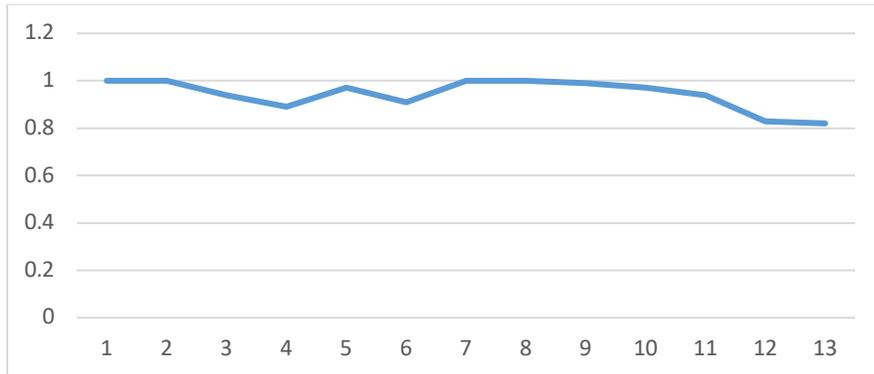

Figure 7. The mean values of DM

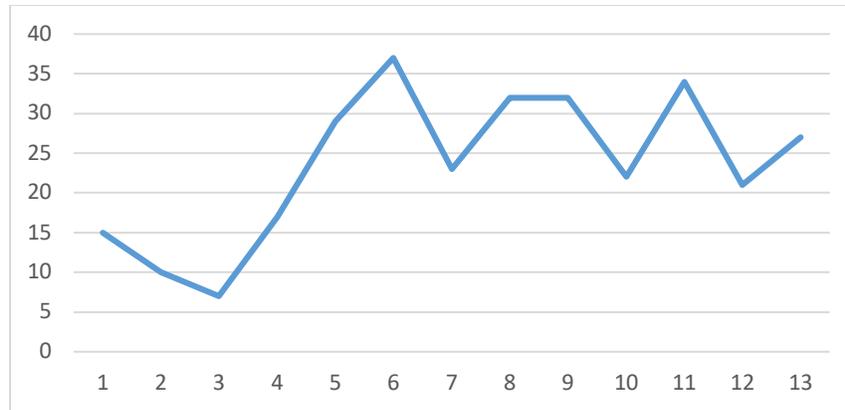
Figure 8. The mean values of NPS

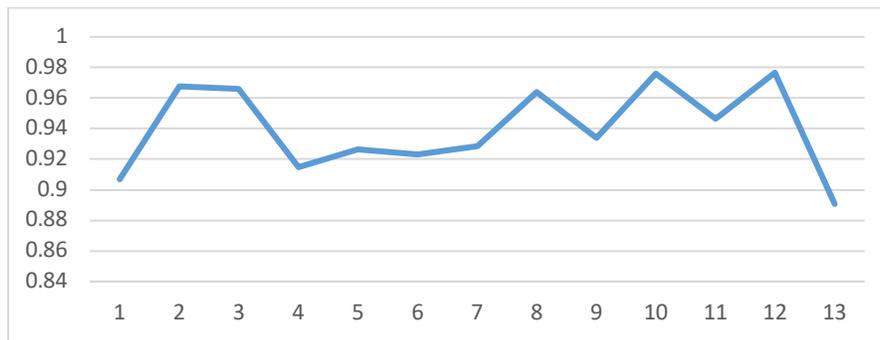
Figure 9. The mean values of MID

**6.3 Solving a numerical example and sensitivity analysis**

Recognizing the inherent potential of carpooling, we conduct a sensitivity analysis concerning the quantity of available vehicles, aiming to discern the optimal number conducive to maximizing the profitability of carpooling. To achieve this, the I9 problem, encompassing 62 nodes, is resolved for varying participant counts. As depicted in the accompanying illustrations, heightened individual participation not only markedly influences organizational expenditures but also efficaciously underscores the beneficial impact of carpooling in curtailing carbon emissions and augmenting incentives.

Within Figure 10, the sensitivity analysis pertaining to the primary objective function vis-à-vis the augmentation of carpooling is scrutinized. It is evident that the objective function's value remains invariant with the presence of up to three vehicles. However, upon surpassing this threshold, there is a pronounced cost diminution, characterized by a precipitous decline. It is postulated that as vehicle numbers swell in more expansive scenarios, this cost reduction will plateau, culminating in a stable trajectory.

Figure 11 delves into the secondary objective function, emphasizing the amplification of carpooling utilization. Meanwhile, Figure 12 elucidates the escalation in pollution attributable to increased carpooling. Notably, while there is an uptick in greenhouse gas emissions concomitant with heightened carpooling adoption, juxtaposing this scenario against one where every employee commutes via private vehicles unequivocally demonstrates that carpooling proliferation invariably results in diminished pollution.

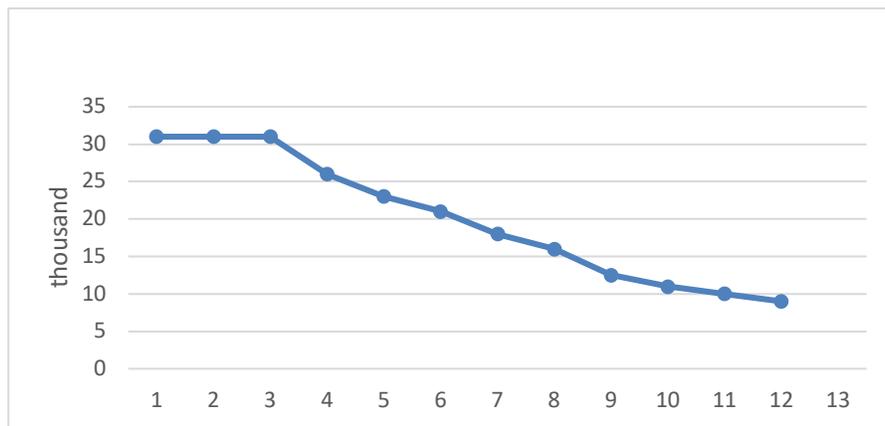

Figure 10. Sensitivity of the first objective function

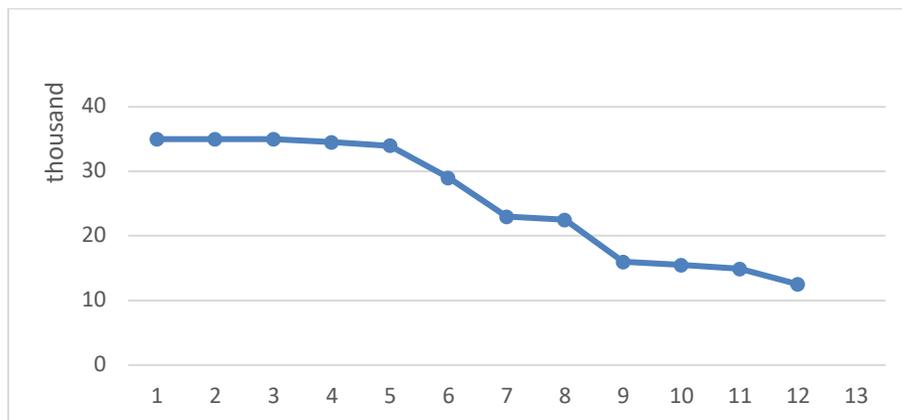

Figure 11. Sensitivity of the second objective function

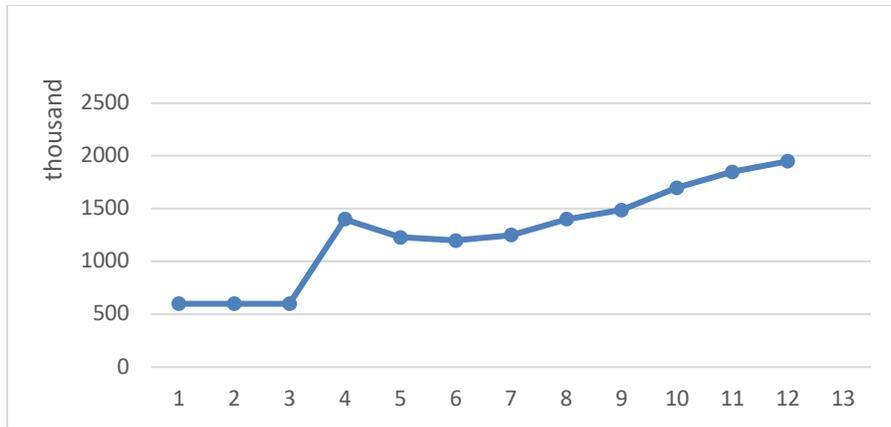
Figure 12. Sensitivity of the third objective function

Owing to the technological strides made in recent epochs, the consumption rates of batteries for hybrid and electric vehicles are anticipated to witness a decline. To elucidate this, subsequent figures delve into this parameter to address the query: how does the diminution in battery consumption rate influence our problem? Figure 13 illustrates that a contraction in battery consumption invariably leads to cost reductions.

Figure 14 conveys that a decrement in the battery consumption rate remains inconsequential to employee satisfaction. Essentially, from the employees' perspective, the quantum of battery consumption lacks significance. Figure 15 indicates that a reduction in the battery consumption rate diminishes the third objective function. This can be attributed to the fact that hybrid and electric vehicles, with reduced battery consumption, possess the capability to traverse extended distances on battery power alone.

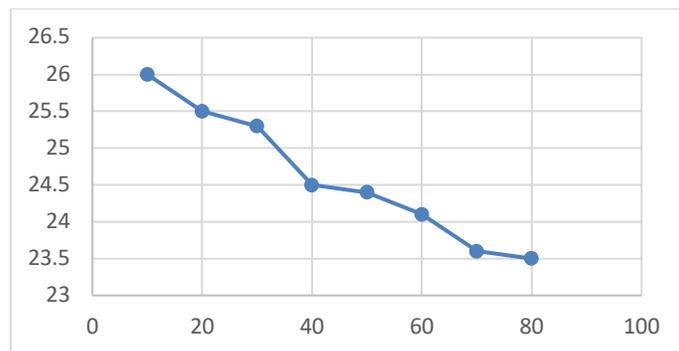
Figure 13. Sensitivity of the first objective function

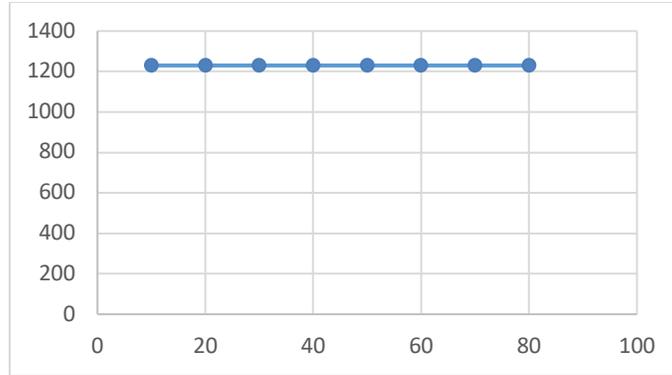
Figure 14. Sensitivity of the second objective function

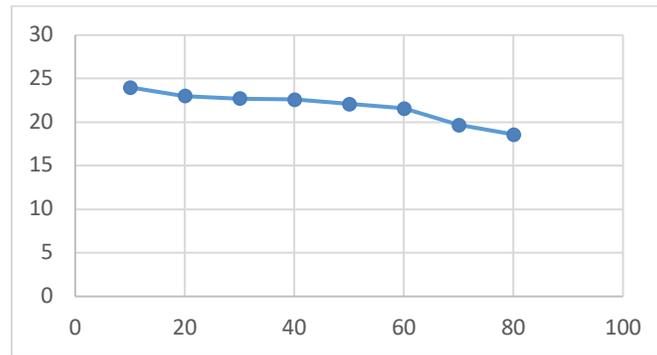
Figure 15. Sensitivity of the third objective function

Computational experiments show that the proposed model can be successfully used by organizations to offer a carpooling transport service to minimize transportation costs. This leads to more employees' involvement and at the same time help raise growing concerns about greenhouse gas emissions by eliminating direct transportation by private cars from employees' homes to the company and vice versa.

## 7. Conclusion and managerial insights

To achieve the goals of this study, we proposed a mathematical model. We used the epsilon constraint method to solve it for small-scale instances and used the NSGA-II algorithm for small-to-large-scale cases. In our comparative experiments, since the results of these two proposed solution methods were close to each other, we concluded that this meta-heuristic algorithm is a suitable for solving our model. The results of numerical experiments showed that carpooling has a significant impact on reducing costs and air pollution compared to the use of private cars and also increases employees' satisfaction.

One research direction is to consider round-trip services for employees. In this study, only one way trips are considered for employees. Another point that can be considered as a future development of this research is that the locations of bus services are not known in advance and suitable locations are located for employees so that their distance to the bus station is minimized. In addition, considering gasoline stations for the fuel-based vehicles and battery charging stations for the electric and hybrid vehicles can be a further development.

**References**

Alumur, S., & Kara, B. Y. (2008). Network hub location problems: The state of the art. European Journal of Operational Research, 190(1), 1-21.

Armbrust, P., Hungerländer, P., Maier, K., & Pachatz, V. (2022). Case study of Dial-a-Ride Problems arising in Austrian rural regions. Transportation Research Procedia, 62, 197-204

Asghari, M., & Al-e, S. M. J. M. (2021). Green vehicle routing problem: A state-of-the-art review. International Journal of Production Economics, 231, 107899.

Asghari, M., Al-e, S. M. J. M., & Rekik, Y. (2022). Environmental and social implications of incorporating carpooling service on a customized bus system. Computers & Operations Research, 142, 105724.

Aziz, M., Oda, T., & Ito, M. (2016). Battery-assisted charging system for simultaneous charging of electric vehicles. Energy, 100, 82-90.

Beraldi, P., De Maio, A., Laganà, D., & Violi, A. (2021). A pick-up and delivery problem for logistics e-marketplace services. *Optimization Letters*, *15*, 1565-1577

Chan, C. C., Wong, Y. S., Bouscayrol, A., & Chen, K. (2009). Powering sustainable mobility: roadmaps of electric, hybrid, and fuel cell vehicles [point of view]. Proceedings of the IEEE, 97(4), 603-607.

Christofides, N., Mingozzi, A., & Toth, P. (1981). Exact algorithms for the vehicle routing problem, based on spanning tree and shortest path relaxations. Mathematical Programming, 20(1), 255-282.

Clarke, G., & Wright, J. W. (1964). Scheduling of vehicles from a central depot to a number of delivery points. Operations research, 12(4), 568-581.

Dantzig, G. B., & Ramser, J. H. (1959). The truck dispatching problem. Management Science, 6(1), 80-91.

Dai, R., Ding, C., Gao, J., Wu, X., & Yu, B. (2022). Optimization and evaluation for autonomous taxi ride-sharing schedule and depot location from the perspective of energy consumption. Applied Energy, 308, 118388

Fisher, M. (1995). Vehicle routing. Handbooks in operations research and management science, 8, 1-33.

Furuhata, M., Dessouky, M., Ordóñez, F., Brunet, M. E., Wang, X., & Koenig, S. (2013). Ridesharing: The state-of-the-art and future directions. Transportation Research Part B: Methodological, 57, 28-46.

Golden, B. L., Wasil, E. A., Kelly, J. P., & Chao, I. (1998). The impact of metaheuristics on solving the vehicle routing problem: algorithms, problem sets, and computational results. Fleet management and logistics (pp. 33-56). Springer, Boston, MA.

Guo, J., Long, J., Xu, X., Yu, M., & Yuan, K. (2022). The vehicle routing problem of intercity ride-sharing between two cities. Transportation Research Part B: Methodological, 158, 113-139

Johnsen, L. C., & Meisel, F. (2022). Interrelated trips in the rural dial-a-ride problem with autonomous vehicles. European Journal of Operational Research, 303(1), 201-219.

Karak, A., & Abdelghany, K. (2019). The hybrid vehicle-drone routing problem for pick-up and delivery services. *Transportation Research Part C: Emerging Technologies*, *102*, 427-449

Laporte, G. (1992). The vehicle routing problem: An overview of exact and approximate algorithms. European journal of operational research, 59(3), 345-358.

Laporte, G., Mercure, H., & Nobert, Y. (1986). An exact algorithm for the asymmetrical capacitated vehicle routing problem. Networks, 16(1), 33-46.

Limmer, S. (2023). Bilevel large neighborhood search for the electric autonomous dial-a-ride problem. Transportation Research Interdisciplinary Perspectives, 21, 100876

Lotfi, S., & Abdelghany, K. (2022). Ride matching and vehicle routing for on-demand mobility services. Journal of Heuristics, 28(3), 235-258

Magnanti, T. L. (1981). Combinatorial optimization and vehicle fleet planning: Perspectives and prospects. Networks, 11(2), 179-213.

Mancini, S. (2017). The hybrid vehicle routing problem. Transportation Research Part C: Emerging Technologies, 78, 1-12.

Mokhtarzadeh, M., Tavakkoli-Moghaddam, R., Vahedi-Nouri, B., & Farsi, A. (2020). Scheduling of human-robot collaboration in assembly of printed circuit boards: a constraint programming approach. International Journal of Computer Integrated Manufacturing, 33(5), 460-473.

Moumou, M., Allaoui, R., & Rhofir, K. (2017, April). Kohonen map approach for vehicle routing problem with pick-up and delivering. In *2017 International Colloquium on Logistics and Supply Chain Management (LOGISTIQUA)* (pp. 183-187). IEEE

Mortazavi, A., Ghasri, M., Haghshenas, H., & Ray, T. (2022a). Adaptive Logit Models for Constructing Feasible Solution for the Dial-a-Ride Problem. Available at SSRN. .

Mortazavi, A., Ghasri, M., & Ray, T. (2022b). A linearly decreasing deterministic annealing algorithm for the multi-vehicle dial-a-ride problem. Plos one, 19(2), e0292683. .

Pouls, M., Ahuja, N., Glock, K., & Meyer, A. (2022). Adaptive forecast-driven repositioning for dynamic ride-sharing. Annals of Operations Research, 1-34.

Rabbani, M., Mokhtarzadeh, M., Manavizadeh, N., & Farsi, A. (2021). Solving a bi-objective mixed-model assembly-line sequencing using metaheuristic algorithms considering ergonomic factors, customer behavior, and periodic maintenance. OPSEARCH, 58(3), 513-539.

Rezgui, D., Siala, J. C., Aggoune-Mtalaa, W., & Bouziri, H. (2019). Application of a variable neighborhood search algorithm to a fleet size and mix vehicle routing problem with electric modular vehicles. Computers & Industrial Engineering, 130, 537-550.

Santos, A., McGuckin, N., Nakamoto, H. Y., Gray, D., & Liss, S. (2011). Summary of travel trends: 2009 national household travel survey (No. FHWA-PL-11-022). United States. Federal Highway Administration.

Schenekemberg, C. M., Chaves, A. A., Coelho, L. C., Guimarães, T. A., & Avelino, G. G. (2022). The dial-a-ride problem with private fleet and common carrier. Computers & Operations Research, 147, 105933

Shahin, A., Hanafizadeh, P., & Hladík, M. (2016). Sensitivity analysis of linear programming in the presence of correlation among right-hand side parameters or objective function coefficients. Central European Journal of Operations Research, 24(3), 563-593.

Souza, A. L., Bernardo, M., Penna, P. H., Pannek, J., & Souza, M. J. (2022). Bi-objective optimization model for the heterogeneous dynamic dial-a-ride problem with no rejects. Optimization Letters, 16(1), 355-374

Stiglic, M., Agatz, N., Savelsbergh, M., & Gradisar, M. (2015). The benefits of meeting points in ride-sharing systems. Transportation Research Part B: Methodological, 82, 36-53.

Su, Y., Dupin, N., & Puchinger, J. (2023). A deterministic annealing local search for the electric autonomous dial-a-ride problem. European Journal of Operational Research, 309(3), 1091-1111.

Tirado, G., & Hvattum, L. M. (2017). Improved solutions to dynamic and stochastic maritime pick-up and delivery problems using local search. *Annals of Operations Research*, *253*, 825-843.

Toth, P., & Vigo, D. (1997). An exact algorithm for the vehicle routing problem with backhauls. Transportation science, 31(4), 372-385.

Toth, P., & Vigo, D. (2002). Models, relaxations and exact approaches for the capacitated vehicle routing problem. Discrete Applied Mathematics, 123(1-3), 487-512.

Vincent, F. Y., Redi, A. P., Hidayat, Y. A., & Wibowo, O. J. (2017). A simulated annealing heuristic for the hybrid vehicle routing problem. Applied Soft Computing, 53, 119-132

Wang, D., Wang, Q., Yin, Y., & Cheng, T. C. E. (2023). Optimization of ride-sharing with passenger transfer via deep reinforcement learning. Transportation Research Part E: Logistics and Transportation Review, 172, 103080

Weiss, M., Zerfass, A., & Helmers, E. (2019). Fully electric and plug-in hybrid cars-An analysis of learning rates, user costs, and costs for mitigating $CO_2$ and air pollutant emissions. Journal of cleaner production, 212, 1478-1489.

Zachariadis, E. E., Tarantilis, C. D., & Kiranoudis, C. T. (2015). The load-dependent vehicle routing problem and its pick-up and delivery extension. *Transportation Research Part B: Methodological*, *71*, 158-181

Zhang, F., Yang, Z. J., Wang, Y., & Kuang, F. (2016). An augmented estimation of distribution algorithm for multi-carpooling problem with time window. In 2016 IEEE 83rd Vehicular Technology Conference (VTC Spring) (pp. 1-5). IEEE.

Kiaghadi, M., & Hoseinpour, P. (2023). University admission process: a prescriptive analytics approach. *Artificial Intelligence Review*, *56*(1), 233-256.

Kiaghadi, M., Sheikholeslami, M., Alinia, A. M., & Boora, F. M. (2024). Predicting the performance of a photovoltaic unit via machine learning methods in the existence of finned thermal storage unit. *Journal of Energy Storage*, *90*, 111766.37